\documentclass[12pt,intlimits,reqno]{article}

\usepackage{amssymb}
\usepackage{epsfig}
\usepackage{graphicx}
\usepackage{lscape}
\usepackage{amsmath}
\usepackage{amsfonts}
\usepackage{amscd}
\usepackage{cite}
\usepackage{amsthm}

\usepackage[all,cmtip]{xy}

\newenvironment{theo}[1]{\vskip+0.4cm {\fontsize{14}{16pt}\selectfont \textbf{Theorem #1.\;}}}{\vskip+0.4cm}

\begin{document}

\title{\textbf{Diophantine Equations and \\  Congruent  Number \\ Equation Solutions}}
\author{\textbf{Mamuka Meskhishvili}}
\date{}
\maketitle

\begin{abstract}
\vskip+0.4cm
By using pairs of nontrivial rational solutions of congruent number equation
\vskip+0.02cm
{\fontsize{14}{16pt}\selectfont
$$  C_N:\;\;y^2=x^3-N^2x,        $$
}
\vskip+0.2cm
\noindent
constructed are pairs of rational right (Pythagorean) triangles with one common side and the other sides equal to the sum and difference of the squares of the same rational numbers.

The parametrizations are found for following Diophantine systems:
\vskip+0.02cm
{\fontsize{14}{16pt}\selectfont
\begin{align*}
    (p^2\pm q^2)^2-a^2 & =\square_{1,2}\,, \\[0.2cm]
    c^2-(p^2\pm q^2)^2 & =\square_{1,2}\,, \\[0.2cm]
    a^2+(p^2\pm q^2)^2 & =\square_{1,2}\,, \\[0.2cm]
    (p^2\pm q^2)^2-a^2 & =(r^2\pm s^2)^2.
\end{align*}
}
\medskip


\medskip

\textbf{2010 AMS Classification.}
    11D25, 11D41, 11D72, 14G05, 14H52.
\end{abstract}

\vskip+2cm

\section{Introduction}
\vskip+0.5cm

By using rational solution {\fontsize{14}{16pt}\selectfont $(x,y)$} of congruent number equation
\vskip+0.02cm
{\fontsize{14}{16pt}\selectfont
$$  C_N:\;\; y^2=x^3-N^2x,      $$
 }
\vskip+0.2cm
\noindent
it is possible to construct a rational right triangle {\fontsize{14}{16pt}\selectfont $(a,b,c)$} with area {\fontsize{14}{16pt}\selectfont $N$} {\fontsize{14}{16pt}\selectfont [1]}
\vskip+0.02cm
{\fontsize{14}{16pt}\selectfont
\begin{equation}\label{eq:1}
\begin{aligned}
    a & =\sqrt{x+N}-\sqrt{x-N}\,, \\[0.2cm]
    b & =\sqrt{x+N}+\sqrt{x-N}\,, \\[0.2cm]
    c & =2\sqrt{x}\,.
\end{aligned}
\end{equation}
}
\vskip+0.02cm

This construction is not possible from every solution. From \eqref{eq:1} follows: a rational right triangle exists, if and only if solution {\fontsize{14}{16pt}\selectfont $x$} satisfies:
\vskip+0.02cm
{\fontsize{14}{16pt}\selectfont
\begin{equation}\label{eq:2}
    x=\square\,, \quad x+N=\square\,, \quad x-N=\square\,.
\end{equation}
}
\vskip+0.2cm
\noindent
These conditions we call Fibonacci conditions, because Fibonacci was the first to find such a number:
\vskip+0.02cm
{\fontsize{14}{16pt}\selectfont
$$  x=\Big(\frac{41}{12}\Big)^2, \;\; N=5.      $$
 }
\vskip+0.2cm
\noindent
Congruent equation solutions which satisfy \eqref{eq:2} we call Fibonacci solutions, respectively. Fibonacci solution has three properties:
\vskip+0.02cm
{\fontsize{14}{16pt}\selectfont
$$  x=\Big(\frac{L}{K}\Big)^2, \;\; K=2M, \;\; (L,N)=1,     $$
 }
\vskip+0.2cm
\noindent
where {\fontsize{14}{16pt}\selectfont $L$}, {\fontsize{14}{16pt}\selectfont $K$} and {\fontsize{14}{16pt}\selectfont $M$} positive integers (as well {\fontsize{14}{16pt}\selectfont $N$}).

If there exists even one nontrivial {\fontsize{14}{16pt}\selectfont $(y\neq 0)$} rational solution of {\fontsize{14}{16pt}\selectfont $C_N$}, that means there exists infinitely many rational solutions. Among them are infinitely many Fibonacci and non-Fibonacci solutions. The tangent and the secant methods can be used to construct infinitely many solutions from the original. Draw the tangent line to {\fontsize{14}{16pt}\selectfont $C_N$} curve at original point {\fontsize{14}{16pt}\selectfont $(x,y)$}, this line will meet the curve in a second point {\fontsize{14}{16pt}\selectfont [1]}:
\vskip+0.02cm
{\fontsize{14}{16pt}\selectfont
$$  \bigg(\Big(\frac{x^2+N^2}{2y}\Big)^2;\frac{(x^2+N^2)(x^4+N^4-6x^2N^2)}{8y^3}\bigg).       $$
 }
\vskip+0.2cm
\noindent
Now using a secant is possible to produce a new rational point and so forth.

In this paper is shown the usage of congruent equation solutions to solve Diophantine systems, including Fibonacci and non-Fibonacci solutions as well. Through them are constructed pairs of rational right (Pythagorean) triangles with one common side and the other sides equal to the sum and difference of the squares of the same rational numbers {\fontsize{14}{16pt}\selectfont $p^2\pm q^2$}.

\vskip+1cm
\section{Common Legs, Hypotenuses $\boldsymbol{p^2\pm q^2}$}
\vskip+0.5cm

The shared (common) legs denote {\fontsize{14}{16pt}\selectfont $a=a_1=a_2$} and hypotheses {\fontsize{14}{16pt}\selectfont $c_{1,2}=p^2\pm q^2$}. This case is equivalent to solving Diophantine system:
\vskip+0.02cm
{\fontsize{14}{16pt}\selectfont
\begin{equation}\label{eq:3}
    \begin{cases}
        (p^2+q^2)^2-a^2=\square\,, \\[0.3cm]
        (p^2-q^2)^2-a^2=\square\,.
    \end{cases}
\end{equation}
}
\vskip+0.2cm
\noindent
We try to find nontrivial integer solutions of system \eqref{eq:3}. Nontrivial solutions mean both Pythagorean triangles exist and they are distinct:
\vskip+0.02cm
{\fontsize{14}{16pt}\selectfont
$$  a_1b_1c_1a_2b_2c_2\neq 0,   $$
}
\vskip+0.2cm
\noindent
or in system notations:
\vskip+0.02cm
{\fontsize{14}{16pt}\selectfont
$$  apq\neq 0, \quad p\neq q.       $$
}
\vskip0.2cm
\noindent
By using rational parametrization formulas for unit hyperbola:
\vskip+0.02cm
{\fontsize{14}{16pt}\selectfont
\begin{equation}\label{eq:4}
    \frac{p^2+q^2}{a}=\frac{1+\xi^2}{2\xi}\,, \quad \frac{p^2-q^2}{a}=\frac{1+\zeta^2}{2\zeta}\,,
\end{equation}
}
\vskip+0.2cm
\noindent
where {\fontsize{14}{16pt}\selectfont $\xi$} and {\fontsize{14}{16pt}\selectfont $\zeta$} are arbitrary nontrivial rationals. After simplification~\eqref{eq:4}
\vskip+0.02cm
{\fontsize{14}{16pt}\selectfont
\begin{equation}\label{eq:5}
\begin{aligned}
    p^2 & =\frac{a(\xi+\zeta)(1+\xi\zeta)}{4\xi\zeta}\,, \\[0.3cm]
    q^2 & =\frac{a(\zeta-\xi)(1-\xi\zeta)}{4\xi\zeta}\,.
\end{aligned}
\end{equation}
}
\vskip+0.2cm
\noindent
From \eqref{eq:5}:
\vskip+0.02cm
{\fontsize{14}{16pt}\selectfont
\begin{align*}
    (\zeta^2-\xi^2)\big(1-(\xi\zeta)^2\big) & =\square\,, \\[0.2cm]
    \Big(\frac{\xi}{\zeta}\Big)(\xi\zeta)\Big(1-\Big(\frac{\xi}{\zeta}\Big)^2\Big)\big(1-(\xi\zeta)^2\big) & =\square\,. \\[0.2cm]
\end{align*}
}
\noindent
By using Property 2 in {\fontsize{14}{16pt}\selectfont [2]} it is possible to express {\fontsize{14}{16pt}\selectfont $\xi$} and {\fontsize{14}{16pt}\selectfont $\zeta$} by arbitrary nontrivial different rational solutions {\fontsize{14}{16pt}\selectfont $x$} and {\fontsize{14}{16pt}\selectfont $z$} {\fontsize{14}{16pt}\selectfont $(xz=\square)$} of arbitrary {\fontsize{14}{16pt}\selectfont $C_N$} congruent number equation:
\vskip+0.02cm
{\fontsize{14}{16pt}\selectfont
$$  \frac{\xi}{\zeta}=\frac{x}{N}\,, \quad \xi\zeta=\frac{z}{N}\,.     $$
 }
\vskip+0.2cm
\noindent
Accordingly,
\vskip+0.02cm
{\fontsize{14}{16pt}\selectfont
\begin{equation}\label{eq:6}
    \xi=\frac{\sqrt{xz}}{N}\,, \quad \zeta=\sqrt{\frac{z}{x}}\,.
\end{equation}
}
\vskip+0.2cm
\noindent
Substitute \eqref{eq:6} in \eqref{eq:5}
\vskip+0.02cm
{\fontsize{14}{16pt}\selectfont
\begin{align*}
    p^2 & =\frac{a}{4N\sqrt{xz}}\,(x+N)(z+N), \\[0.4cm]
    q^2 & =\frac{a}{4N\sqrt{xz}}\,(x-N)(z-N).
\end{align*}
}

\begin{theo}{1}
\textit{Nontrivial integer solutions of system \eqref{eq:3} are given by formulas:
\vskip+0.02cm
{\fontsize{14}{16pt}\selectfont
\begin{align*}
    a & =k\cdot 4N\sqrt{xz}\,, \\[0.2cm]
    p^2 & =k\cdot (x+N)(z+N), \\[0.2cm]
    q^2 & =k\cdot (x-N)(z-N),
\end{align*}
}
\vskip+0.2cm
\noindent
where {\fontsize{14}{16pt}\selectfont $k$} is integer, {\fontsize{14}{16pt}\selectfont $x$} and {\fontsize{14}{16pt}\selectfont $z$} are nontrivial different rational solutions of arbitrary congruent number equation.}
\end{theo}

Integer {\fontsize{14}{16pt}\selectfont $k$} and solutions {\fontsize{14}{16pt}\selectfont $x$}, {\fontsize{14}{16pt}\selectfont $z$} must be chosen so that {\fontsize{14}{16pt}\selectfont $a$}, {\fontsize{14}{16pt}\selectfont $p$}\, and {\fontsize{14}{16pt}\selectfont $q$} are integers.

Direct calculations give:
\vskip+0.02cm
{\fontsize{14}{16pt}\selectfont
\begin{equation}\label{eq:7}
\begin{aligned}
    a_1 & =k\cdot 4N\sqrt{xz}\,, \\[0.2cm]
     b_1 & = k\cdot 2(xz-N^2), \\[0.2cm]
    c_1 & =p^2+q^2=k\cdot 2(xz+N^2); \\[0.2cm]
    a_2 & =k\cdot 4N\sqrt{xz}\,, \\[0.2cm]
     b_2 & = k\cdot 2N(x-z), \\[0.2cm]
    c_2 & =p^2-q^2=k\cdot 2N(x+z).
\end{aligned}
\end{equation}
}

\newpage
\subsection*{Numerical Examples}

To construct examples of Diophantine system \eqref{eq:3}, we use {\fontsize{14}{16pt}\selectfont $C_N$} solution from~{\fontsize{14}{16pt}\selectfont [1]}.

\begin{enumerate}
\item[1)] {\fontsize{14}{16pt}\selectfont $N=6$};\quad {\fontsize{14}{16pt}\selectfont $x=18$}, {\fontsize{14}{16pt}\selectfont $z=\dfrac{19602}{47^2}$}\,.
\vskip+0.02cm
{\fontsize{14}{16pt}\selectfont
\begin{align*}
    (74^2+23^2)^2-4653^2 & =3796^2, \\[0.2cm]
    (74^2-23^2)^2-4653^2 & =1680^2.
\end{align*}
}
\vskip+0.2cm

\item[2)] {\fontsize{14}{16pt}\selectfont $N=34$};\quad {\fontsize{14}{16pt}\selectfont $x=162$}, {\fontsize{14}{16pt}\selectfont $z=\dfrac{2178}{7^2}$}\,.
\vskip+0.02cm
{\fontsize{14}{16pt}\selectfont
\begin{align*}
    (217^2+64^2)^2-35343^2 & =37024^2, \\[0.2cm]
    (217^2-64^2)^2-35343^2 & =24480^2.
\end{align*}
}
\end{enumerate}

\vskip+1cm
\section{Common Hypotheses (Siblings), Legs $\boldsymbol{p^2\pm q^2}$}
\vskip+0.5cm

The hypotenuses denote {\fontsize{14}{16pt}\selectfont $c=c_1=c_2$} and legs {\fontsize{14}{16pt}\selectfont $a_{1,2}=p^2\pm q^2$}. Diophantine system is:
\vskip+0.02cm
{\fontsize{14}{16pt}\selectfont
\begin{equation}\label{eq:8}
    \begin{cases}
        c^2-(p^2+q^2)^2=\square\,, \\[0.3cm]
        c^2-(p^2-q^2)^2=\square\,.
    \end{cases}
\end{equation}
}
\vskip+0.2cm
\noindent
Again we seek nontrivial integer solutions of system \eqref{eq:8}
\vskip+0.02cm
{\fontsize{14}{16pt}\selectfont
\begin{gather*}
    a_1b_1c_1a_2b_2c_2\neq 0; \\[0.2cm]
    cpq\neq 0, \quad p\neq q.
\end{gather*}
}
\vskip+0.2cm
\noindent
By using parametrization formulas for unit circle:
\vskip+0.02cm
{\fontsize{14}{16pt}\selectfont
\begin{equation}\label{eq:9}
    \frac{p^2+q^2}{c}=\frac{2\xi}{1+\xi^2}\,, \quad \frac{p^2-q^2}{a}=\frac{2\zeta}{1+\zeta^2}\,.
\end{equation}
}
\vskip+0.2cm
\noindent
After simplification \eqref{eq:9}:
\vskip+0.02cm
{\fontsize{14}{16pt}\selectfont
\begin{equation}\label{eq:10}
\begin{aligned}
    p^2 & =\frac{c(\xi+\zeta)(1+\xi\zeta)}{(1+\xi^2)(1+\zeta^2)}\,, \\[0.4cm]
    q^2 & =\frac{c(\xi-\zeta)(1-\xi\zeta)}{(1+\xi^2)(1+\zeta^2)}\,.
\end{aligned}
\end{equation}
}
\vskip+0.2cm
\noindent
From \eqref{eq:10}
\vskip+0.02cm
{\fontsize{14}{16pt}\selectfont
$$  (\xi^2-\zeta^2)\big(1-(\xi\zeta)^2\big)=\square\,.        $$
 }
\vskip+0.2cm
\noindent
Again using Property 2 in {\fontsize{14}{16pt}\selectfont [2]}:
\vskip+0.02cm
{\fontsize{14}{16pt}\selectfont
\begin{equation}\label{eq:11}
    \xi=\frac{\sqrt{xz}}{N}\,, \quad \zeta=\sqrt{\frac{z}{x}}\,,
\end{equation}
}
\vskip+0.2cm
\noindent
where {\fontsize{14}{16pt}\selectfont $x$} and {\fontsize{14}{16pt}\selectfont $z$} are nontrivial different rational solutions of arbitrary congruent number equation.

Substitute \eqref{eq:11} in \eqref{eq:10}:
\vskip+0.02cm
{\fontsize{14}{16pt}\selectfont
\begin{align*}
    p^2 & =\frac{c\sqrt{xz}\,(x+N)(z+N)}{(x+z)(xz+N^2)}\,, \\[0.4cm]
    q^2 & =\frac{c\sqrt{xz}\,(x-N)(z-N)}{(x+z)(xz+N^2)}\,.
\end{align*}
}
\vskip+0.2cm

So,

\begin{theo}{2}
\textit{Nontrivial integer solutions of system \eqref{eq:8} are given by formulas:
\vskip+0.02cm
{\fontsize{14}{16pt}\selectfont
\begin{align*}
    c & =k\cdot\sqrt{xz}\,(x+z)(xz+N^2), \\[0.2cm]
    p^2 & =k\cdot xz(x+N)(z+N), \\[0.2cm]
    q^2 & =k\cdot xz(x-N)(z-N),
\end{align*}
}
\vskip+0.2cm
\noindent
where {\fontsize{14}{16pt}\selectfont $k$} is integer, {\fontsize{14}{16pt}\selectfont $x$} and {\fontsize{14}{16pt}\selectfont $z$} are nontrivial different rational solutions of arbitrary congruent number equation.}
\end{theo}

Integer {\fontsize{14}{16pt}\selectfont $k$} and solutions {\fontsize{14}{16pt}\selectfont $x$}, {\fontsize{14}{16pt}\selectfont $z$} must be chosen so that {\fontsize{14}{16pt}\selectfont $c$}, {\fontsize{14}{16pt}\selectfont $p$}\, and {\fontsize{14}{16pt}\selectfont $q$} are integers.

Direct calculations give:
\vskip+0.02cm
{\fontsize{14}{16pt}\selectfont
\begin{equation}\label{eq:12}
\begin{aligned}
    a_1 & =p^2+q^2=k\cdot 2xz(xz+N^2), \\[0.2cm]
    b_1 & =k\cdot \sqrt{xz}\,(x-z)(xz+N^2), \\[0.2cm]
    c_1 & =k\cdot \sqrt{xz}\,(x+z)(xz+N^2); \\[0.2cm]
    a_2 & =p^2-q^2=k\cdot 2Nxz(x+z), \\[0.2cm]
    b_2 & =k\cdot \sqrt{xz}\,(x+z)(xz-N^2), \\[0.2cm]
    c_2 & =k\cdot \sqrt{xz}\,(x+z)(xz+N^2).
\end{aligned}
\end{equation}
}

\newpage
\subsection*{Numerical Examples}

To construct examples of Diophantine system \eqref{eq:8}, we use the same solutions.

\begin{enumerate}
\item[1)] {\fontsize{14}{16pt}\selectfont $N=6$};\quad {\fontsize{14}{16pt}\selectfont $x=18$}, {\fontsize{14}{16pt}\selectfont $z=\dfrac{19602}{47^2}$}\,.
\vskip+0.02cm
{\fontsize{14}{16pt}\selectfont
\begin{align*}
    15358381995^2-(114774^2+35673^2)^2 & =5215702800^2, \\[0.2cm]
    15358381995^2-(114774^2-35673^2)^2 & =9708645804^2.
\end{align*}
}
\vskip+0.2cm

\item[2)] {\fontsize{14}{16pt}\selectfont $N=34$};\quad {\fontsize{14}{16pt}\selectfont $x=162$}, {\fontsize{14}{16pt}\selectfont $z=\dfrac{2178}{7^2}$}\,.
\vskip+0.02cm
{\fontsize{14}{16pt}\selectfont
\begin{align*}
    3322469535^2-(50127^2+14784^2)^2 & =1891797600^2, \\[0.2cm]
    3322469535^2-(50127^2-14784^2)^2 & =2403264864^2.
\end{align*}
}
\end{enumerate}

\vskip+1cm
\section{Common Legs, Other Legs $\boldsymbol{p^2\pm q^2}$}
\vskip+0.5cm

Denote common legs {\fontsize{14}{16pt}\selectfont $a=a_1=a_2$} and other legs {\fontsize{14}{16pt}\selectfont $b_{1,2}=p^2\pm q^2$}. Diophantine system is
\vskip+0.02cm
{\fontsize{14}{16pt}\selectfont
\begin{equation}\label{eq:13}
    \begin{cases}
        a^2+(p^2+q^2)^2=\square\,, \\[0.3cm]
        a^2+(p^2-q^2)^2=\square\,.
    \end{cases}
\end{equation}
}
\vskip+0.2cm
\noindent
We seek nontrivial integer solutions of system \eqref{eq:13}:
\vskip+0.02cm
{\fontsize{14}{16pt}\selectfont
\begin{gather*}
    a_1b_1c_1a_2b_2c_2\neq 0; \\[0.2cm]
    apq\neq 0, \quad p\neq q.
\end{gather*}
}
\vskip+0.2cm

Rational parametrization for unit hyperbola gives:
\vskip+0.02cm
{\fontsize{14}{16pt}\selectfont
\begin{equation}\label{eq:14}
    \frac{p^2+q^2}{a}=\frac{1-\xi^2}{2\xi}\,, \quad \frac{p^2-q^2}{a}=\frac{1-\zeta^2}{2\zeta}\,.
\end{equation}
}
\vskip+0.2cm
\noindent
Accordingly from \eqref{eq:14}:
\vskip+0.02cm
{\fontsize{14}{16pt}\selectfont
\begin{equation}\label{eq:15}
\begin{aligned}
    p^2 & =\frac{a(\xi+\zeta)(1-\xi\zeta)}{4\xi\zeta}\,, \\[0.4cm]
    q^2 & =\frac{a(\zeta-\xi)(1+\xi\zeta)}{4\xi\zeta}\,.
\end{aligned}
\end{equation}
}
\vskip+0.2cm
\noindent
From \eqref{eq:15}
\vskip+0.02cm
{\fontsize{14}{16pt}\selectfont
$$  (\zeta^2-\xi^2)\big(1-(\xi\zeta)^2\big)=\square\,.        $$
 }
\vskip+0.2cm
\noindent
Again using Property 2 in {\fontsize{14}{16pt}\selectfont [2]}:
\vskip+0.02cm
{\fontsize{14}{16pt}\selectfont
$$  \xi=\frac{\sqrt{xz}}{N}\,, \quad \zeta=\sqrt{\frac{z}{x}}\,,        $$
 }
\vskip+0.2cm
\noindent
and after substitution in \eqref{eq:15}
\vskip+0.02cm
{\fontsize{14}{16pt}\selectfont
\begin{align*}
    p^2 & =\frac{a(x+N)(z-N)}{4N\sqrt{xz}}\,, \\[0.4cm]
    q^2 & =\frac{a(x-N)(z+N)}{4N\sqrt{xz}}\,.
\end{align*}
}

\newpage
\begin{theo}{3}
\textit{Nontrivial integer solutions of system \eqref{eq:13} are given by formulas:
\vskip+0.02cm
{\fontsize{14}{16pt}\selectfont
\begin{align*}
    a & =k\cdot4N\sqrt{xz}\,, \\[0.2cm]
    p^2 & =k\cdot (x+N)(z-N), \\[0.2cm]
    q^2 & =k\cdot (x-N)(z+N),
\end{align*}
}
\vskip+0.2cm
\noindent
where {\fontsize{14}{16pt}\selectfont $k$} is integer, {\fontsize{14}{16pt}\selectfont $x$} and {\fontsize{14}{16pt}\selectfont $z$} are nontrivial different rational solutions of arbitrary congruent number equation.}
\end{theo}

Integer {\fontsize{14}{16pt}\selectfont $k$} and solutions {\fontsize{14}{16pt}\selectfont $x$}, {\fontsize{14}{16pt}\selectfont $z$} must be chosen so that {\fontsize{14}{16pt}\selectfont $a$}, {\fontsize{14}{16pt}\selectfont $p$}\, and {\fontsize{14}{16pt}\selectfont $q$} are integers.

Direct calculations give:
\vskip+0.02cm
{\fontsize{14}{16pt}\selectfont
\begin{equation}\label{eq:16}
\begin{aligned}
    a_1 & =k\cdot 4N\sqrt{xz}\,, \\[0.2cm]
    b_1 & =p^2+q^2=k\cdot 2(xz-N^2), \\[0.2cm]
    c_1 & =k\cdot 2(xz+N^2); \\[0.2cm]
    a_2 & =k\cdot 4N\sqrt{xz}\,, \\[0.2cm]
    b_2 & =p^2-q^2=k\cdot 2N(x-z), \\[0.2cm]
    c_2 & =k\cdot 2N(x+z).
\end{aligned}
\end{equation}
}
\vskip+0.2cm

Solutions of system \eqref{eq:13} are impossible to find using the preceding congruent number equation solutions, so we have to take other {\fontsize{14}{16pt}\selectfont $C_N$} solutions.

\newpage
\subsection*{Numerical Examples}
\vskip+0.2cm

\begin{enumerate}
\item[1)] {\fontsize{14}{16pt}\selectfont $N=5$};\quad {\fontsize{14}{16pt}\selectfont $x=\dfrac{12005}{31^2}$}\,, {\fontsize{14}{16pt}\selectfont $z=45$}.
\vskip+0.02cm
{\fontsize{14}{16pt}\selectfont
\begin{align*}
    4557^2+(82^2+60^2)^2 & =11285^2, \\[0.2cm]
    4557^2+(82^2-60^2)^2 & =5525^2.
\end{align*}
}
\vskip+0.2cm

\item[2)] {\fontsize{14}{16pt}\selectfont $N=34$};\quad {\fontsize{14}{16pt}\selectfont $x=\dfrac{833}{4^2}$}\,, {\fontsize{14}{16pt}\selectfont $z=\dfrac{153}{2^2}$}\,.
\vskip+0.02cm
{\fontsize{14}{16pt}\selectfont
\begin{align*}
    1344^2+(17^2+9^2)^2 & =1394^2, \\[0.2cm]
    1344^2+(17^2-9^2)^2 & =1360^2.
\end{align*}
}
\end{enumerate}

\vskip+1cm
\section{Common Legs, Hypotenuses $\boldsymbol{p^2\pm q^2}$, \\ Other Legs $\boldsymbol{r^2\pm s^2}$}
\vskip+0.5cm

Preceding numerical examples are constructed by using non-Fibonacci solutions. Of course if we use Fibonacci solution {\fontsize{14}{16pt}\selectfont $(x,z)$} it is possible to construct solutions for each Diophantine systems \eqref{eq:3}, \eqref{eq:8} and \eqref{eq:13}. It indicates that intersection case may exist.

Consider Pythagorean triangles with sides:
\vskip+0.02cm
{\fontsize{14}{16pt}\selectfont
$$  (a,r^2+s^2,p^2+q^2), \quad (a,r^2-s^2,p^2-q^2),      $$
 }
\vskip+0.2cm
\noindent
accordingly
\vskip+0.02cm
{\fontsize{14}{16pt}\selectfont
$$  a=a_1=a_2, \quad b_{1,2}=r^2\pm s^2, \quad c_{1,2}=p^2\pm q^2.      $$
 }
\vskip+0.2cm
\noindent
Corresponding Diophantine system is:
\vskip+0.02cm
{\fontsize{14}{16pt}\selectfont
\begin{equation}\label{eq:17}
    \begin{cases}
        (p^2+q^2)^2-a^2=(r^2+s^2)^2, \\[0.3cm]
        (p^2-q^2)^2-a^2=(r^2-s^2)^2;
    \end{cases}
\end{equation}
}
\vskip+0.2cm
\noindent
again we seek nontrivial integer solutions:
\vskip+0.02cm
{\fontsize{14}{16pt}\selectfont
$$  apqrs\neq 0, \quad p\neq q, \quad r\neq s.      $$
}
\vskip+0.2cm
\noindent
Divide system \eqref{eq:17} by two subsystems:
\vskip+0.02cm
{\fontsize{14}{16pt}\selectfont
$$  \begin{cases}
        (p^2+q^2)^2-a^2=\square\,, \\[0.3cm]
        (p^2-q^2)^2-a^2=\square
    \end{cases} \quad \text{and} \quad
        \begin{cases}
            a^2+(r^2+s^2)^2=\square\,, \\[0.3cm]
            a^2+(r^2-s^2)^2=\square\,.
        \end{cases}         $$
}
\vskip+0.2cm
\noindent
For these subsystems use parametric formulas from Theorem 1 and Theorem~3, respectively:
\vskip+0.02cm
{\fontsize{14}{16pt}\selectfont
\begin{equation}\label{eq:18}
\begin{aligned}
    a & =k_1\cdot 4N_1\sqrt{x_1z_1}\,, \\[0.2cm]
    p^2 & =k_1\cdot (x_1+N_1)(z_1+N_1), \\[0.2cm]
    q^2 & =k_1\cdot (x_1-N_1)(z_1-N_1); \\[0.2cm]
    a & =k_2\cdot 4N_2\sqrt{x_2z_2}\,, \\[0.2cm]
    r^2 & =k_2\cdot (x_2+N_2)(z_2-N_2), \\[0.2cm]
    s^2 & =k_2\cdot (x_2-N_2)(z_2+N_2).
\end{aligned}
\end{equation}
}
\vskip+0.2cm

Expressions \eqref{eq:7} and \eqref{eq:16} gives:
\vskip+0.02cm
{\fontsize{14}{16pt}\selectfont
\begin{align}
    k_1\cdot 4N_1\sqrt{x_1z_1} & =k_2\cdot 4N_2\sqrt{x_2z_2}, \tag{$19_1$} \label{eq:19_1} \\[0.2cm]
    k_1\cdot 2(x_1z_1-N_1^2) & =k_2\cdot 2(x_2z_2-N_2^2), \tag{$19_2$} \label{eq:19_2} \\[0.2cm]
    k_1\cdot 2(x_1z_1+N_1^2) & =k_2\cdot 2(x_2z_2+N_2^2), \tag{$19_3$} \label{eq:19_3} \\[0.2cm]
    k_1\cdot 2N_1(x_1-z_1) & =k_2\cdot 2N_2(x_2-z_2), \tag{$19_4$} \label{eq:19_4} \\[0.2cm]
    k_1\cdot 2N_1(x_1+z_1) & =k_2\cdot 2N_2(x_2+z_2). \tag{$19_5$} \label{eq:19_5}
\end{align}
}
\vskip+0.2cm
\noindent
The sum and difference \eqref{eq:19_2} and \eqref{eq:19_3}, \eqref{eq:19_4} and \eqref{eq:19_5} gives:
\vskip+0.02cm
{\fontsize{14}{16pt}\selectfont
\begin{align}
    k_1N_1^2 & =k_2N_2^2, \tag{$20_1$} \label{eq:20_1} \\[0.2cm]
    k_1x_1z_1 & =k_2x_2z_2, \tag{$20_2$} \label{eq:20_2} \\[0.2cm]
    k_1N_1x_1 & =k_2N_2x_2, \tag{$20_3$} \label{eq:20_3} \\[0.2cm]
    k_1N_1z_1 & =k_2N_2z_2. \tag{$20_4$} \label{eq:20_4}
\end{align}
\begin{align}
    \text{\eqref{eq:20_3}, \eqref{eq:20_1}} &\;\; \Longrightarrow \;\;\frac{x_1}{x_2}=\frac{k_2N_2}{k_1N_1}=\frac{N_1}{N_2}\,,
                \tag{$21_1$} \label{eq:21_1} \\[0.4cm]
    \text{\eqref{eq:20_4}, \eqref{eq:20_1}} &\;\; \Longrightarrow \;\;\frac{z_1}{z_2}=\frac{k_2N_2}{k_1N_1}=\frac{N_1}{N_2}\,.
                \tag{$21_2$} \label{eq:21_2}
\end{align}
}
\vskip+0.2cm

From \eqref{eq:21_1} follows, exists rational {\fontsize{14}{16pt}\selectfont $t$}
\vskip+0.02cm
{\fontsize{14}{16pt}\selectfont
$$  N_1=x_1t \quad\text{and}\quad N_2=x_2t.     $$
 }
\vskip+0.2cm
\noindent
Because {\fontsize{14}{16pt}\selectfont $x_1$} and {\fontsize{14}{16pt}\selectfont $x_2$} are solutions {\fontsize{14}{16pt}\selectfont $C_{N_1}$} and {\fontsize{14}{16pt}\selectfont $C_{N_2}$}, respectively:
\setcounter{equation}{21}
\vskip+0.02cm
{\fontsize{14}{16pt}\selectfont
\begin{equation}\label{eq:22}
\begin{aligned}
    x_1(x_1^2-N_1^2)=x_1^3(1-t^2) & =\square\,, \\[0.2cm]
    x_2(x_2^2-N_2^2)=x_2^3(1-t^2) & =\square\,.
\end{aligned}
\end{equation}
}
\vskip+0.2cm
\noindent
From \eqref{eq:22} follows
\vskip+0.02cm
{\fontsize{14}{16pt}\selectfont
\begin{equation}\label{eq:23}
    \frac{x_1}{x_2}=\square\,.
\end{equation}
}
\vskip+0.2cm
\noindent
Expressions \eqref{eq:23} and \eqref{eq:21_1} gives
\vskip+0.02cm
{\fontsize{14}{16pt}\selectfont
$$  \frac{N_1}{N_2}=\square\,.        $$
 }
\vskip+0.2cm
\noindent
Because congruent numbers are squafree so,
\vskip+0.02cm
{\fontsize{14}{16pt}\selectfont
\begin{equation}\label{eq:24}
    N_1=N_2.
\end{equation}
}
\vskip+0.2cm
\noindent
From \eqref{eq:21_1}, \eqref{eq:21_2} and \eqref{eq:20_1} follows
\vskip+0.02cm
{\fontsize{14}{16pt}\selectfont
\begin{equation}\label{eq:25}
    x_1=x_2, \quad z_1=z_2, \quad k_1=k_2.
\end{equation}
}
\vskip+0.2cm
\noindent
Conditions \eqref{eq:24} and \eqref{eq:25} give parametric solution of this intersection case:

\begin{theo}{4}
\textit{Nontrivial integer solutions of system \eqref{eq:17} are given by formulas:
\vskip+0.02cm
{\fontsize{14}{16pt}\selectfont
\begin{equation}\label{eq:26}
\begin{aligned}
    p^2 & =k\cdot (x+N)(z+N), \\[0.2cm]
    q^2 & =k\cdot (x-N)(z-N), \\[0.2cm]
    a & =k\cdot4N\sqrt{xz}\,, \\[0.2cm]
    r^2 & =k\cdot (x+N)(z-N), \\[0.2cm]
    s^2 & =k\cdot (x-N)(z+N),
\end{aligned}
\end{equation}
}
\vskip+0.2cm
\noindent
where {\fontsize{14}{16pt}\selectfont $k$} is integer, {\fontsize{14}{16pt}\selectfont $x$} and {\fontsize{14}{16pt}\selectfont $z$} are nontrivial different rational solutions of arbitrary congruent number equation.}
\end{theo}

Integer {\fontsize{14}{16pt}\selectfont $k$} and solutions {\fontsize{14}{16pt}\selectfont $x$}, {\fontsize{14}{16pt}\selectfont $z$} must be chosen so that {\fontsize{14}{16pt}\selectfont $p$}, {\fontsize{14}{16pt}\selectfont $q$}, {\fontsize{14}{16pt}\selectfont $a$}, {\fontsize{14}{16pt}\selectfont $r$}, {\fontsize{14}{16pt}\selectfont $s$} are integers. Because Fibonacci solutions satisfy \eqref{eq:2}, from them it is always possible to construct integer solution of \eqref{eq:17}.

Direct calculations give:
\vskip+0.02cm
{\fontsize{14}{16pt}\selectfont
\begin{align*}
    a_1 & =k\cdot 4N\sqrt{xz}\,, \\[0.2cm]
    b_1 & =k\cdot 2(xz-N^2), \\[0.2cm]
    c_1 & =k\cdot 2(xz+N^2); \\[0.2cm]
    a_2 & =k\cdot 4N\sqrt{xz}\,, \\[0.2cm]
    b_2 & =k\cdot 2N(x-z), \\[0.2cm]
    c_2 & =k\cdot 2N(x+z).
\end{align*}
}

\subsection*{Numerical Examples}

\begin{enumerate}
\item[1)] Non-Fibonacci pair: \\[0.2cm]
     {\fontsize{14}{16pt}\selectfont $N=5$};\quad {\fontsize{14}{16pt}\selectfont $x=\Big(\dfrac{5}{2}\Big)^2$}, {\fontsize{14}{16pt}\selectfont$z=\Big(\dfrac{41}{12}\Big)^2$}.
\vskip+0.02cm
{\fontsize{14}{16pt}\selectfont
\begin{align*}
    (735^2+155^2)^2-492000^2 & =(465^2+245^2)^2, \\[0.2cm]
    (735^2-155^2)^2-492000^2 & =(465^2-245^2)^2.
\end{align*}
}
\vskip+0.2cm

\item[2)] Fibonacci pair:\\[0.2cm]
    {\fontsize{14}{16pt}\selectfont $N=6$};\quad {\fontsize{14}{16pt}\selectfont $x=\Big(\dfrac{5}{2}\Big)^2$}, {\fontsize{14}{16pt}\selectfont $z=\Big(\dfrac{1201}{140}\Big)^2$}.
\vskip+0.02cm
{\fontsize{14}{16pt}\selectfont
\begin{align*}
    (8743^2+1151^2)^2-40353600^2 & =(8057^2+1249^2)^2, \\[0.2cm]
    (8743^2-1151^2)^2-40353600^2 & =(8057^2-1249^2)^2.
\end{align*}
}
\end{enumerate}

Diophantine system \eqref{eq:17} is equivalent to:
\vskip+0.02cm
{\fontsize{14}{16pt}\selectfont
\begin{equation}\label{eq:27}
    \begin{cases}
        p^4+q^4-a^2=r^4+s^4, \\[0.2cm]
        pq=rs.
    \end{cases}
\end{equation}
}
\vskip+0.2cm
\noindent
By Theorem 4 the second equation is satisfied automatically. This means: the parametric formulas in \eqref{eq:26} are solution of the first equation in system~\eqref{eq:27}.

\begin{theo}{5}
\textit{Nontrivial integer solutions of equation
\vskip+0.02cm
{\fontsize{14}{16pt}\selectfont
\begin{equation}\label{eq:28}
     p^4+q^4-a^2=r^4+s^2,
\end{equation}
}
\vskip+0.2cm
\noindent
which satisfy \,{\fontsize{14}{16pt}\selectfont $pq=rs$}\, condition, are given by parametric formulas \eqref{eq:26}.}
\end{theo}

Obviously there are solutions \eqref{eq:28} which do not satisfy \,{\fontsize{14}{16pt}\selectfont $pq=rs$}\, condition. For example, it is possible to construct such solutions by using the equal sum of biquadratics  {\fontsize{14}{16pt}\selectfont[3]}:
\vskip+0.02cm
{\fontsize{14}{16pt}\selectfont
\begin{align*}
    7^4+28^4 & =3^4+20^4+26^4, \\[0.2cm]
    51^4+76^4 & =5^4+42^4+78^4, \\[0.2cm]
    37^4+38^4 & =25^4+26^4+42^4.
\end{align*}
}

\subsection*{Numerical Examples}

{\fontsize{14}{16pt}\selectfont
\begin{align*}
    735^4+155^4-492000^2 & =465^4+245^4, \\[0.2cm]
    8743^4+1151^4-40353600^2 & =8057^4+1249^4.
\end{align*}
}

\newpage
\begin{center}
     {\fontsize{14}{16pt}\selectfont \textbf{References} }
\end{center}
\vskip+0.5cm

\begin{enumerate}
\item[1.] Meskhishvili M., Three-Century Problem. \emph{Tbilisi}, 2013.

\medskip

\item[2.] Meskhishvili M., Perfect Cuboid and Congruent Number Equation Solutions. 2013, \emph{Unsolved Problems in Number Theory, Logic, and Cryptography}, \texttt{www.unsolvedproblems.org}; \\ \texttt{arxiv.org/pdf/1211.6548v2.pdf}.

\item[3.] Dickson L. E., History of the Theory of Numbers, Volume II: Diophantine Analysis. \emph{Dover, New York}, 2005.

\end{enumerate}

\vskip+1cm

\noindent \textbf{Author's address:}

\medskip

\noindent {Georgian-American High School, 18 Chkondideli Str., Tbilisi 0180, Georgia.}

\noindent {\small \textit{E-mail:} \texttt{director@gahs.edu.ge} }

\end{document}